\documentclass[leqno,11pt]{amsart}
\usepackage{amssymb,amsmath,latexsym,amsfonts,amsbsy, amsthm,esint}
\usepackage{color}
\usepackage{graphicx}
\usepackage{tikz}
\usepackage{tikz-3dplot}
\usetikzlibrary{arrows, calc, decorations.markings,intersections}
\usepackage{caption}
\usepackage{subcaption}
\usepackage{hyperref,cleveref}
\usepackage[foot]{amsaddr}
\usepackage[normalem]{ulem}

\let\pa=\partial
\let\al=\alpha

\let\de=\delta
\let\z=\zeta
\let\lam=\lambda

\let\f=\frac

\let\Om=\Omega

\let\e=\varepsilon

\let\ri=\rightarrow
\let\na=\nabla

\newtheorem{theorem}{Theorem}[section]
\newtheorem{lemma}{Lemma}[section]
\newtheorem{prop}{Proposition}[section]

\theoremstyle{remark}

\newtheorem{case}{Case}
\newtheorem{rmk}{Remark}[section]

\newcommand{\dist}{\mathrm{dist}}

\newcommand{\BS}{\mathbb{S}^2}

\newcommand{\bn}{\mathbf{n}}
\newcommand{\bx}{\mathbf{x}}
\newcommand{\bp}{\mathbf{p}}
\newcommand{\bq}{\mathbf{Q}}

\newcommand{\bm}{\mathbf{m}}

\begin{document}
\title{Eigenframe discontinuities of the Q-tensor model}
\author{Zhiyuan Geng$^1$}
\email{geng42@purdue.edu}
\author{Changyou Wang$^1$}
\email{wang2482@purdue.edu}
\address{$^1$Department of Mathematics, Purdue University, 150 N. University Street, West Lafayette, IN 47907}

\date{\today}

\begin{abstract}   
In this paper, we study the defect structure of minimizer of a Landau-de Gennes energy functional in three-dimensional domains, subject to constraint $|Q|=1$. The set of defects is identified by discontinuities in both the eigenframe and the leading eigenvector. Through a blow-up analysis, we prove that the defect set is 1-rectifiable and classify the asymptotic profile of the leading eigenvector near singularities. This generalizes some previous results on the structure of ring disclinations in the $Q$-tensor model. 

\end{abstract}

\maketitle

\section{Introduction}\label{sec:intro}

\subsection{Mathematical formulation and main results}

Liquid crystal phases are intermediate states of matter between crystalline solids and isotropic fluids, exhibiting characteristics of both. The main types of liquid crystals include nematic, smectic, and cholesteric phases. Among them, nematic liquid crystals (NLC) are the most studied mathematically. They consist of rigid-rod shaped molecules that locally align in the same direction. Sharp variations in this alignment, known as ``defect patterns", are commonly observed in experiments and typically appear as isolated points or disclination lines. Defects are essential in liquid crystal display technologies, as they influence optical properties and can be strategically engineered to control light modulation. Mathematically, they represent singular structures in orientation of molecules and give rise to many challenging problems in analysis and geometry. Readers are referred to survey articles \cite{Ball, Lin-Liu} and references therein for a more detailed discussion.

Defect patterns are predicted and described differently across various continuum models. In the Oseen-Frank theory \cite{frank}, the local average orientation of NLC molecules is represented by a unit-vector field that takes value in $\BS$, with defects defined as discontinuities of the director $\bn$. The structure and size of the defects in Oseen-Frank theory have been widely studied, see \cite{su,bcl,hkl1,hkl2}. A major limitation of the Oseen-Frank model is that it only allows for point defects (in three-dimensional domains) with finite energy.

The Ericksen model \cite{ericksen1991liquid} addresses this limitation by relaxing the unit-length constraint and introducing a scalar order parameter $s\in [-\f12,1]$. In this framework, defects are interpreted as regions of "isotropic melting", where the degree of orientation $s$ vanishes. In certain parameter regime, the defect set consists of a 1-rectifiable curve along which the director field $\bn$ forms a planar half-degree configuration, see \cite{lin1989nonlinear,lin1991nematic,hardt1993harmonic,alper2017defects,alper2018rectifiability}.  The Ericksen model allows for both point and line defects. 

In the Landau-de Gennes theory \cite{deGennes}, the order parameter is represented by a symmetric traceless $3\times 3$ matrix called \emph{$Q$-tensor}, which takes values in the space
\begin{equation*}
    \mathcal{Q}_0:=\{Q\in \mathbb{R}^{3\times 3}:\  Q=Q^T,\, \mathrm{tr}(Q)=0\}.
\end{equation*}
$\mathcal{Q}_0$ can be identified with the Euclidean space $\mathbb{R}^5$.

When $Q=0$, it is called an isotropic phase. It is uniaxial if two of the nonzero eigenvalues are identical and it is biaxial if $Q$ has three distinct eigenvalues. Biaxiality can be measured by a signed biaxiality parameter $\beta$ \cite{dipasquale1}, which is defined by 
\begin{equation}\label{def:beta}
    \beta(Q):=\sqrt{6}\f{\mbox{tr}(Q^3)}{|Q|^3}\in[-1,1],\quad \text{for }Q\neq 0.
\end{equation}
For a nonzero $Q$, when $\beta(Q)=1$, $Q$ has two equal negative eigenvalues and one positive eigenvalue. This configuration is called \emph{positive uniaxial}, indicating that most liquid crystal polymers align along the direction of the leading eigenvector, i.e. the eigenvector corresponding to the only positive eigenvalue. On the other hand, if $\beta(Q)=-1$, $Q$ is called \emph{negative uniaxial}, characterized by two equal positive eigenvalues and one negative eigenvalue. Near the negative uniaxial point, the liquid crystal polymers tend to align perpendicularly to the eigenvector associated with the negative eigenvalue. In the case where $\beta(Q)\in (-1,1)$, $Q$ is biaxial with three distinct eigenvalues. 

Let $\Omega\subset \mathbb{R}^3$ be a simply connected domain occupied by nematic liquid crystal materials, the simplest form of the Landau-de Gennes free energy for $Q:\Omega\rightarrow\mathcal{Q}_0$ is given by
\begin{equation}
    \label{ene:LdG}
    \mathcal{E}_{LdG}(Q):=\int_{\Omega} \left(\f{L}{2}|\na Q|^2+f_{b}(Q)\right)\,d\bx,
   \ \   Q\in H^1(\Omega,\mathcal{Q}_0)
\end{equation}
where the bulk potential $f_b(Q)$ takes the form
\begin{equation}
    \label{ene:bulk potential}
    f_b(Q):=-\f{a^2}{2}\mathrm{tr}(Q^2)-\f{b^2}{3}\mathrm{tr}(Q^3)+\f{c^2}{4}[\mathrm{tr}(Q^2)]^2+C.
\end{equation}
Here $L,a,b,c$ are material-dependent constants and $C=C(a,b,c)$ is a constant ensuring $\min\limits_{Q\in \mathcal{Q}_0} f_b(Q)=0$. It is well-known that $f_b$ attains its minimum value on a submanifold of $\mathcal{Q}_0$ defined by
\begin{equation}
    \label{def: vacuum}
    \mathcal{N}:=\{Q\in\mathcal{Q}_0:\, Q=s_+(\bn\otimes \bn-\f13\mathrm{Id}),\ \bn\in\mathbb{S}^2\},
\end{equation}
where $\mathrm{Id}$ is the $3\times 3$ identity matrix and
\begin{equation*}
    s_+:=\f{b^2+\sqrt{b^4+24a^2c^2}}{4c^2}.
\end{equation*}
Note that $\beta(Q)=1$ for $Q\in \mathcal{N}$. A minimizer $Q$ of the Landau-de Gennes energy 
$\mathcal{E}_{LdG}$ in $H^1(\Omega,\mathcal{Q}_0)$, subject to the Dirichlet boundary condition 
\begin{equation}
    \label{Dirichlet bdy cond}
    Q|_{\pa\Om}=Q_b(\bx)\in C^\infty(\pa \Om, \mathcal{Q}_0), 
\end{equation}
satisfies the Euler-Lagrange equation
\begin{equation}\label{ELQ}
    L\Delta Q=-a^2Q-b^2[Q^2-\f13\mathrm{tr}(Q)^2\,\mathrm{Id}]+c^2\mathrm{tr}(Q)^2Q.
\end{equation}

As the target manifold $\mathcal{Q}_0$ is isometric to $\mathbb{R}^5$, 
any minimizer $Q$ of $\mathcal{E}_{LdG}$ is a real analytic function. Therefore, the defect sets are often interpreted as the regions where \emph{eigenframe discontinuities} or \emph{eigenvalue crossings} occur, see e.g. \cite{de1972types, palffy1994new,sonnet1995alignment}. The precise definition of the defect set in our setting is given by \eqref{def:defect} in Section \ref{sec:preliminary}.

We rewrite the Landau-de Gennes energy in the same form as in \cite{dipasquale1}. Rescale the Q-tensor by setting 
\begin{equation*}
    \mathbf{Q}(\bx):=\sqrt{\f{3}{2}}\f{1}{s_+} Q({\bf x}),
\end{equation*}
then the energy functional rewrites
\begin{equation*}
    \mathcal{E}_{LdG}(Q)=(\f23s_+^2L)\cdot \mathcal{E}_{\lambda,\mu}(\mathbf{Q}),
\end{equation*}
with
\begin{equation}\label{def: E_{lambda,mu}}
    \mathcal{E}_{\lambda,\mu}(\mathbf{Q}):= \int_{\Om} \left( \f12|\na \mathbf{Q}|^2+\lambda W(\mathbf{Q})+\f{\mu}{4}(1-|\bq|^2)^2 \right)\,d\bx.
\end{equation}
Here the parameters $\lambda$ and $\mu$ are defined by 
\begin{equation*}
    \lambda:=\sqrt{\f23}\f{b^2s_+}{L},\quad \mu:=\f{a^2}{L},
\end{equation*}
and $W(\mathbf{Q})$ takes the form
\begin{equation}\label{def: W(Q)}
    W(\bq):=\f{1}{4\sqrt{6}}|\bq|^4-\f13\mathrm{tr}(\bq^3)+\f{1}{12\sqrt{6}}.
\end{equation}
In particular, when $|\bq|=1$ i.e. $\bq\in \mathbb{S}^4$, $W(\bq)$ becomes
\begin{equation}\label{W(Q)when|Q|=1}
    W(\bq)=\f{1}{3\sqrt{6}} (1-\beta(\bq)),
\end{equation}
which clearly only penalizes biaxiality. 
For simplicity, we write $\mathbf{Q}$ as $Q$ and assume $\Omega=B_1$, the unit ball in $\mathbb{R}^3$, from now on.

Following \cite{lyu,dipasquale1}, by setting $\lambda=1$ and $\mu\ri\infty$ we consider the minimizer $Q\in H^{1}(B_1,\mathbb{S}^4)$ of the following energy functional
\begin{equation}\label{ene}
    \mathcal{E}(Q):=\int_{B_1}\left(\f12|\na Q|^2+\f{1}{3\sqrt{6}}(1-\beta(Q))\right)\,d\bx,
\end{equation}
subject to the Dirichlet boundary condition
\begin{equation}\label{bdy cond}
\begin{split}
    &\qquad Q|_{\pa B_1} =Q_b\in C^\omega(\pa B_1;\mathbb{S}^4),\\
    &\min\limits_{\bx\in \pa B_1}\beta(Q(\bx))\geq -1+\delta\ \ \text{ for some }\delta>0,\quad \mathrm{deg}(\bn(\bx),\pa B_1)\neq 0.
    \end{split}
\end{equation}
Here $\bn(\bx)$ denotes the unit eigenvector associated with the largest eigenvalue $\lambda_1$ of $Q(\bx)$. On $\pa B_1$, since $\beta(Q(\bx))>-1$, the largest eigenvalue $\lambda_1$ of $Q(\bx)$ is simple so that $\bn(\bx)$ is a smooth vector field on $\pa B_1$. The constraint $|Q|=1$ is called the ``Lyuksyutov constraint" in the series of works by Dipasquale et al. \cite{dipasquale1,dipasquale2,dipasquale3}, which physically corresponds to the zero temperature limit. We also point out that all results can extend naturally to any simply connected domain $\Omega$ without difficulty.  

We summarize our main result in the following theorem, whose proofs will be given by Proposition \ref{prop: k=1 profile}, \ref{prop: rectifiable S1} and \ref{prop: U, k>=2}.

\begin{theorem}
    Let $Q\in H^1(B_1,\mathbb{S}^4)$ be a minimizer of the energy functional $\mathcal{E}(Q)$ given by \eqref{ene}, which satisfies the boundary condition \eqref{bdy cond}. The defect set $\mathcal{S}(Q)$, where both the eigenframe and the leading eigenvector are discontinuous (see \eqref{def:defect}), is countably 1-rectifiable. For $\mathcal{H}^1$ almost everywhere on $\mathcal{S}(Q)$, the blow-up profile of $Q$ is a function of two variables after an appropriate rotation. In addition, they can be classified according to the vanishing order of the function $Q+\f{\sqrt{6}}{2}(\bp(\bx)\otimes \bp(\bx)-\f13\mathrm{Id})$, where $\bp(\bx)$ denotes the eigenvector associated with the negative eigenvalue.  
\end{theorem}

Our results provide a description of the set where ``eigenvalue crossing" occurs and the leading eigenvector fails to be smooth in the framework of $Q$-tensor model with the Lyuksyutov constraint. This generalizes the ``ring defect" configurations in \cite{yu2020,ty}, where half-degree disclination lines are obtained under the assumption of axial symmetry. 

In the proof, we essentially consider a minimizing harmonic map $Q:B_1\ri \mathbb{S}^4$, noting that the potential energy $W(Q)$ is a lower order term and will vanish after blow-ups. We analyze the blow-up profile near the point where the eigenvalues of $Q$ are $\f{\sqrt{6}}{6},\f{\sqrt{6}}{6},-\f{\sqrt{6}}{3}$. After subtracting the leading-order term $-\f{\sqrt{6}}{2}(\bp(\bx)\otimes \bp(\bx)-\f13\mathrm{Id})$ from $Q$, the remaining function characterizes the asymptotic behavior of the leading eigenvector in the blow-up limit. We show that any tangent map is a ``harmonic function" with a source term determined by $\bp(\bx)$. This enables us to adapt techniques from the study of singular sets for elliptic equations \cite{hs,han1994singular,han1998geometric,hardt1999critical} to classify the tangent maps.

\subsection{Related results}

In recent years, defect configurations in Landau-de Gennes model have been extensively studied through mathematical analysis. In the so-called large domain regime, the energy-minimizing configurations are investigated when $L\rightarrow 0$, or equivalently, $\lambda$ and $\mu$ tend to $\infty$ at the same rate. This corresponds to the case where the defect core size becomes vanishingly small, compared to the domain size. As $L\ri 0$, the tensor field $Q_L$ strongly converges to a minimizing harmonic map $Q_0: \Omega\ri \mathcal{N}$ outside a small neighborhood of the defect set $\mathcal{S}$. When the total energy $\mathcal{E}_{LdG}(Q_L)$ is uniformly bounded, $\mathcal{S}$ consists of a set of  finite points in dimension three \cite{mz,nz,gz}; while when $\mathcal{E}_{LdG}(Q_L)$ scales as $O(|\ln{L}|)$, $\mathcal{S}$ is a set of co-dimension two, i.e. a set of finite points in dimension two and a 1-rectifiable set in 
dimension three \cite{bpp,gm,canevari1,canevari2}. The techniques in these works are similar to those applied in the study of the Ginzburg-Landau model \cite{bbh93,bbh94,LW1,LW2}.

A complete characterization of the profile inside the defect core remains largely open. Many studies focus on typical defect core configurations and their stability. Among these, the melting hedgehog solution, given by  
\begin{equation*}
    Q(\bx) = r(\bx) \left( \frac{\bx}{|\bx|} \otimes \frac{\bx}{|\bx|} - \frac{1}{3} I_d \right),
\end{equation*}  
has been widely investigated \cite{ss, gm1999, ma2012, lamy2013, insz1, insz2}. This configuration is uniaxial everywhere and vanishes at the origin. However, in certain parameter regimes, the radial hedgehog becomes unstable, necessitating the emergence of biaxiality near the defect core. This phenomenon, known as ``biaxial escape," has been rigorously analyzed in the low temperature regime, that is, when $\mu$ is much larger than $\lambda$ \cite{insz2, cl2017, hmp}.  To provide an intuitive explanation, we examine the energy functional \eqref{def: E_{lambda,mu}}. 
The term $\lambda W(Q) $ penalizes biaxiality, whereas the term $ \frac{\mu}{4} (1 - |Q|^2)^2 $ penalizes deviations of $Q$ from $ \mathbb{S}^4 $. When $ \mu $ is significantly larger than $\lambda$, the potential energy penalizes less biaxiality, making the biaxial state energetically more favorable than isotropic melting. For a more detailed discussion of the mechanisms behind biaxial escape, see \cite{gartland2018scalings}.

There are two primary types of biaxial core structure: the half-degree ring disclination and the split-core solution. These biaxial configurations were first discovered and studied numerically \cite{mg2000,hqz} and have recently been rigorously constructed in \cite{yu2020,dipasquale1,ty,dipasquale2,dipasquale3} in the axially symmetric setting. In a series of articles by Dipasquale et al. \cite{dipasquale1,dipasquale2,dipasquale3}, the authors consider minimizers of $\mathcal{E}(Q)$ under the Lyuksyutov constraint and the assumption of axial symmetry. The minimizer $Q$ is smooth except for a possible finite set of singularities on the axis of symmetry. In particular, the set $\{\beta(Q)=-1\}$ contains a circle, and $\{\beta=t\}$ for any $t\in(-1,1)$ is a union of finitely many axially symmetric tori. In \cite{yu2020}, ring disclination and split-core solutions are obtained by considering a constrained minimization problem under more restrictive symmetry hypotheses. The results are further extended from the low-temperture limit (Lyuksyutov regime) to the general Landau-de Gennes model in \cite{ty}. 

Our results are closely related to these aforementioned works. In fact, for the torus solution $Q$ in \cite{dipasquale1,dipasquale2,dipasquale3}, the set $\{\beta(Q)=-1\}$ is only shown to be nonempty for topological reasons, and therefore contains a circle due to the axial symmetry assumption. However, neither the size of the set nor the behavior of the eigenframe near the negative uniaxial sets has been established. In \cite{yu2020,ty}, the detailed configurations of eigenvalues and eigenframes are derived near the ring disclinations where $\beta(Q)=-1$, but these behaviors heavily rely on the extra symmetry assumptions, namely, that one of the eigenvectors is fixed and the solution is equivariant under reflection across $\{z=0\}$. In this paper, we aim to provide a better understanding of the size and eigenframe profiles near negative uniaxial sets $\{\beta(Q)=-1\}$ without imposing any symmetry constraints.

The paper is organized as follows. In Section \ref{sec:preliminary}, we review the regularity properties of the minimizer $Q$ and introduce the order parameter $U\approx Q+\f{\sqrt{6}}{2}(\bp\otimes\bp-\f13\mathrm{Id})$ which encodes information about local behaviors of the eigenframe.  Section \ref{sec:functional} formulates the energy functional and the blow-up sequence for $U$. In Section \ref{sec:classification}, we classify the blow-up profiles for $U$ based on its vanishing order and show that the discontinuity set for the eigenframe is countably 1-rectifiable, thus proving the main theorem.

\section{Preliminaries} \label{sec:preliminary}

Throughout this paper, when we say a quantity $A(s)\sim O(s^k)$ for an integer $k$ and the independent variable $s$, we mean that there exists $c\geq 0$ such that 
\begin{equation*}
 \lim\limits_{s\ri 0+} \f{|A(s)|}{s^k}=c.
\end{equation*}

Define the admissible space for $Q:\Omega\to\mathcal{Q}_0$ 
by
$$
\mathcal{A}_{Q_b}:=\{Q\in H^{1}(\Omega, \mathcal{Q}_0),\  |Q|=1, \ Q=Q_b \text{ on }\pa B_1, \text{ with }Q_b  \text{ satisfying }\eqref{bdy cond}\}.
$$
We first recall the regularity results on a minimizer $Q\in \mathcal{A}_{Q_b}$ of the energy functional $\mathcal{E}(Q)$ defined by \eqref{ene}. The proof can be essentially reduced to the classical Schoen-Uhlenbeck theory on the regularity of harmonic maps from a three-dimensional domain into $\mathbb{S}^4$ \cite{su,su2}. This result was recently reproved in \cite[Theorem 1.2]{dipasquale1}, specifically adapted to the Q-tensor model. We omit the proof and refer interested readers to \cite{su,su2, dipasquale1}.

\begin{prop} Let $Q\in \mathcal{A}_{Q_b}$ be a minimizer of the functional $\mathcal{E}(\cdot)$ defined in \eqref{ene}. Then $Q\in C^\omega(\bar{B}_1)$.
\end{prop}

Consequently, the defect set is identified with the singular structure of the eigenframe of tensor $Q$. For any point $\bx\in B_1$, $Q(\bx)$ has the expression 
\begin{equation}\label{decomp in eigv}
    Q(\bx)=\lam_1(\bx)(\bn(\bx)\otimes \bn(\bx))+\lam_2(\bx)(\bm(\bx)\otimes\bm(\bx))+\lam_3(\bx)(\bp(\bx)\otimes \bp(\bx)),
\end{equation}
where $\lam_1\geq \lam_2\geq \lam_3$ are eigenvalues of $Q(\bx)$ and $\bn(\bx)$, $\bm(\bx)$, $\bp(\bx)$ are corresponding unit eigenvectors. Since $Q\in \mathcal{Q}_0$ satisfies $|Q|=1$, we have 
\begin{equation}\label{cond on eigenvalue}
\sum_{i=1}^3\lam_i(\bx)=0,\quad \sum_{i=1}^3\lam_i^2(\bx)=1, \quad \forall \bx\in\Omega.    
\end{equation}

The biaxiality parameter function $\beta$ defined in \eqref{def:beta} can be rewritten in terms of eigenvalues:
\begin{equation*}
    \beta(Q)=\sqrt{6}\sum\limits_{i=1}^3 \lambda_i^3.
\end{equation*}

It is elementary to show that $\beta(Q)=1$ if and only if $\lambda_1=\f{\sqrt{6}}{3}$ and  $\lambda_2=\lambda_3=-\f{\sqrt{6}}{6}$; while $\beta(Q)=-1$ if and only if $\lambda_1=\lambda_2=\f{\sqrt{6}}{6}$ and $\lambda_3=-\f{\sqrt{6}}{3}$. When $\beta(Q)\in(-1,1)$, $Q$ has three distinct eigenvalues. By the well-known fact that  simple eigenvalues and their associated eigenvectors inherit the smoothness of $Q(\bx)$ (see for instance \cite{nomizu1973characteristic}), we conclude that eigenframe singularities belong to the set $\{\beta=\pm 1\}$. Furthermore, when $\beta(Q(\bx_0))=1$, the leading eigenvector $\bn(\bx)$ remains smooth near $\bx_0$, meaning that the primary preferred direction of nematic liquid crystal polymers does not have a singular behavior near $\bx_0$. Thus we do not identify positive uniaxial points as defects. Instead, our focus is on the negative uniaxial points (i.e. $\beta(Q)=-1$), where the two positive eigenvalues exchange and the leading eigenvector may exhibit discontinuity. This leads to the following definition of the defect set. For $Q\in \mathcal{A}_{Q_b}$, define
\begin{equation}
    \label{def:defect}
    \mathcal{S}(Q):=\Big\{\bx\in B_1: \beta(Q(\bx))=-1,\ \text{the eigenframe }\{\bn,\bm,\bp\} \text{ is discontinuous at }\bx\Big\}.
\end{equation}

Take an arbitrary point $\bx_0\in \mathcal{S}(Q)$, the decomposition \eqref{decomp in eigv} writes
\begin{equation*}
    Q(\bx_0)=\f{\sqrt{6}}{6}(\bn(\bx_0)\otimes \bn(\bx_0))+\f{\sqrt{6}}{6}(\bm(\bx_0)\otimes\bm(\bx_0))-\f{\sqrt{6}}{3}(\bp(\bx_0)\otimes \bp(\bx_0)).
\end{equation*}

We have the following lemma regarding the eigenvalue decomposition of $Q$ in a small neighborhood of $\bx_0$.

\begin{lemma}\label{lemma: eigenvalue of Q}  Given a minimizer $Q\in \mathcal{A}_{Q_b}$ of $\mathcal{E}(Q)$.
For sufficiently small $\e>0$, there exists $r_0=r_0(\e)>0$ such that the three ordered eigenvalues $\{\lam_i(\bx)\}_{i=1}^3$ of $Q(\bx)$ satisfy 
\begin{equation*} 
\begin{split}
& \qquad\qquad  \lam_3(\bx)\in C^{\omega}(B_{r_0}(\bx_0)), \quad \lam_1(\bx),\, \lam_2(\bx) \in C^{0,1}_{loc}(B_{r_0}(\bx_0)).\\ 
&0\leq \lam_1(\bx)-\frac{\sqrt{6}}{6}\leq \e,\; -\e\leq \lam_2(\bx)-\frac{\sqrt{6}}{6}\leq 0,\; 0\leq \lam_3(\bx)+\frac{\sqrt{6}}{3}\leq \e,\quad \forall \bx\in B_{r_0}(\bx_0).
\end{split}
\end{equation*}
Moreover, there exists a positive integer $k$ and a constant $c>0$ such that 
\begin{equation}\label{vanish order of lam_1,lam_3}
    \begin{split}
    c&=\lim\limits_{r\ri 0} \f{\fint_{ B_r(\bx_0)}\left(\lam_1(\bx)-\frac{\sqrt{6}}{6}\right)^2\,d\bx}{r^{2k}},\\
    \f{\sqrt{6}}{3}c&=\lim\limits_{r\ri 0} \f{\fint_{ B_r(\bx_0)}\left(\lam_3(\bx)+\frac{\sqrt{6}}{3}\right)\,d\bx}{|r|^{2k}}.
    \end{split}
\end{equation}
\end{lemma}

\begin{proof}
    From \eqref{cond on eigenvalue} and $\lam_1\geq \lam_2\geq \lam_3$ it is straightforward to deduce that 
    \begin{equation*}
        \lam_3\geq -\f{\sqrt{6}}{3},\quad \lam_1\geq \frac{\sqrt{6}}{6}. 
    \end{equation*}
    As argued before, the isolated eigenvalue $\lam_3(\bx)$ is as smooth as $Q(\bx)$ thanks to the implicit function theorem applied to the characteristic polynomial. Thus, $\lam_3(\bx)$ is analytic near $\bx_0$. The local Lipschitz regularity of $\lam_1(\bx),\,\lam_2(\bx)$ follows from \cite[Theorem 2]{bronshtein1979smoothness}.

    For $\bx$ near $\bx_0$, we set $\delta(\bx):=\lam_3(\bx)+\frac{\sqrt{6}}{3}\geq 0$. Then we derive the perturbations of $\lam_1,\lam_2$ around $\f{\sqrt{6}}{6}$ in terms of $\delta$. Let
    \begin{equation*}
        s(\bx)=\lam_1(\bx)-\f{\sqrt{6}}{6}\geq 0,\quad r(\bx)=\lam_2(\bx)-\f{\sqrt{6}}{6}\le 0.
    \end{equation*}
The relation \eqref{cond on eigenvalue} implies that 
\begin{equation}\label{relation: s,r}
\begin{split}
    &s+r+\delta=0,\\
    &s^2+r^2+\delta^2+\f{\sqrt{6}}{3}(s+r-2\delta)=0.
\end{split}
\end{equation}
A direct computation gives
\begin{equation}\label{solve for s,r}
    s=\left(\f32\right)^{\f14}\sqrt{\de}-\f{\de}{2}+O(\de^{\f32}),\quad r=-\left(\f32\right)^{\f14}\sqrt{\de}-\f{\de}{2}+O(\de^{\f32}),
\end{equation}
which further implies that $\lam_2\leq \frac{\sqrt{6}}{6}$. 

It remains to prove the existence of the integer $k$. Since $\lam_3(\bx)+\f{\sqrt{6}}{3}\geq 0$ is real analytic and vanishes at $\bx_0$, there is an even number $2k$ such that the leading term in the power series expansion of $\lam_3(\bx)+\f{\sqrt{6}}{3}$ near $\bx_0$ is of the order $|\bx-\bx_0|^{2k}$. From \eqref{solve for s,r} we know that 
\begin{equation*}
    \lam_3(\bx)+\f{\sqrt{6}}{3}=\delta\approx\f{\sqrt{6}}{3}s^2=\f{\sqrt{6}}{3}\left(\lam_1(\bx)-\f{\sqrt{6}}{6}\right)^2, 
\end{equation*}
then \eqref{vanish order of lam_1,lam_3} follows immediately. 
\end{proof}

A consequence of Lemma \ref{lemma: eigenvalue of Q} and its proof is the following decomposition of $Q(\bx)$ in a neighborhood $B_r(\bx_0)$:
\begin{equation}\label{formula: Q}
    Q(\bx)=-\f{\sqrt{6}}{2}\left(\bp(\bx)\otimes \bp(\bx)-\f13 \mathrm{Id}\right)+U(\bx)+R(\bx),
\end{equation}
with $U(\bx)$ defined by 
\begin{equation}\label{def:U}
    U(\bx)=s(\bx)(\bn(\bx)\otimes \bn(\bx)-\bm(\bx)\otimes \bm(\bx)).
\end{equation}
Here $\bp(\bx)$ is the eigenvector associated with $\lam_3(\bx)$, $\bn(\bx)$ is the leading eigenvector associated with $\lam_1(\bx)$, $s(\bx)=\f{|U|}{\sqrt{2}}=\lam_1(\bx)-\f{\sqrt{6}}{6}\geq 0$, and $R(\bx)$ is the remainder term of order $O(s^2)$ which by \eqref{solve for s,r} has the form
\begin{equation}\label{formula: R}
 \begin{split}
     R(\bx)&= \left(\lambda_2(\bx)-\frac{\sqrt{6}}{6}+s(\bx)\right)\bm(\bx)\otimes \bm(\bx)+ \left(\lambda_3(\bx)+\frac{\sqrt{6}}{3}\right) \bp(x)\otimes \bp(\bx)\\
     &=\left(\f{\sqrt{6}}{3}s^2+\tau(s)\right)\left(\bp(\bx)\otimes \bp(\bx)-\bm(\bx)\otimes \bm(\bx)\right)\\
     &=\left(\f{\sqrt{6}}{3}s+\f{\tau(s)}{s}\right)\left(\f{U}{2}+\f{3s}{2}(\bp(\bx)\otimes \bp(\bx)-\f13\mathrm{Id})\right),
    \end{split}
\end{equation}
where $\tau(s)\sim O(s^3)$ and can be calculated explicitly through the relation \eqref{relation: s,r}. 

We have that 
\begin{align*}
&\lam_1=\frac{\sqrt{6}}6+s,\ \  \lam_2(\bx)=\f{\sqrt{6}}{6}-s-\f{\sqrt{6}}{3}s^2-\tau(s),\ \  \lam_3(\bx)=-\frac{\sqrt{6}}{3}+\f{\sqrt{6}}{3}s^2+\tau(s),\\
&\qquad \bp(\bx),\, Q(\bx)\in C^\omega(B_r(\bx_0)),\quad  s(\bx)\geq 0,\quad  s(\bx_0)=0.
 \end{align*}
 The vanishing order of the function $U(\bx)$ at $\bx_0$, whose existence is guaranteed by Lemma \ref{lemma: eigenvalue of Q}, is defined as
\begin{equation}\label{def: k}
k(\bx_0):=\f12\lim\limits_{r\ri 0+} \frac{\ln{\fint_{B_r(x_0)}s(x)^2\,dx}}{\ln{r}}.
\end{equation}

We focus on the behavior of $U(\bx)$ near $\bx_0$, treating $\bp(\bx)$ as a given real analytic function. By \eqref{def:U}, $s(\bx)$, $\bn(\bx)$, $\bm(\bx)$ can be recovered from $U(\bx)$ as its eigenvalue and associated eigenvectors. Combining with \eqref{formula: R}, this implies that $R(\bx)$ can be expressed explicitly in terms of $U(\bx)$ and $\bp(\bx)$. Therefore, we can consider $U(\bx)$ as the only order parameter locally. From the definition (\ref{def:U}), $U(\bx)$ only takes value in the space
\begin{equation}\label{target U}
\mathcal{U}_\bp:=\{U\in \mathcal{Q}_0:\,  U(\bx)\bp(\bx)=\mathbf{0}\}.
\end{equation}

{Since $Q(\bx)+\f{\sqrt{6}}{2}\left(\bp(\bx)\otimes \bp(\bx)-\f13\mathrm{Id}\right)$ is a real analytic function in $B_r(\bx_0)$, and $U(\bx)$ contains the leading-order term in its Taylor series expansion at $\bx_0$ by Lemma \ref{lemma: eigenvalue of Q}, \eqref{def:U} and \eqref{formula: R},  it follows that for any $i,j\in \{1,2,3\}$, there exists a degree $k(\bx_0)$ homogeneous polynomial $\bar{U}_{ij}$ such that} 
\begin{equation}\label{polynomial repres}
    U_{ij}(\bx)=\bar{U}_{ij}(\bx-\bx_0)+O(|\bx-\bx_0|^{k(\bx_0)+1}).
\end{equation}

\section{Energy functional for $U$}\label{sec:functional}
We will derive the equation satisfied by $\bar{U}$ by computing the energy functional for which $U$ locally minimizes. Direct computation gives
\begin{equation}\label{formula:gradient Q}
    \na Q(\bx)=-\f{\sqrt{6}}2(\na \bp\otimes \bp+\bp\otimes \na \bp)+\na U(\bx)+\nabla R(\bx).
\end{equation}
\begin{equation}\label{formula:energy density in U}
\begin{split}
     &\f12|\na Q|^2+ \f{1}{3\sqrt{6}}(1-\beta(Q))\\
     =& \f12|\na Q-\na R|^2+(\na Q-\na R):\na R+\f12|\na R|^2+\f{1}{3\sqrt{6}}\big(1-\sqrt{6}(\lam_1^3+\lam_2^3+\lam_3^3)\big)\\
     =&\f32|\na \bp|^2+\f12|\na 
 U|^2-\sqrt{6}\sum\limits_{1\leq i,j,k\leq 3} \pa_i \bp_j\,\bp_k \,\pa_iU_{jk}\\
  &+(\na Q-\na R):\na R+\f12|\na R|^2+ \f{1}{3\sqrt{6}}(2-9s^2)+\zeta(s) \\
  =& \f32|\na \bp|^2+\f12|\na 
 U|^2+\sqrt{6}\sum\limits_{1\leq i,j,k\leq 3} \pa_i \bp_j\,\pa_i \bp_k \,U_{jk}\\
  &+(\na Q-\na R):\na R+\f12|\na R|^2+ \f{1}{3\sqrt{6}}(2-\f{9}{2}|U|^2)+\zeta(s) 
\end{split}
\end{equation}
where $\zeta(s)\sim O(s^3)$. Here we have used the following relations in the above calculation.
\begin{equation*}
0=\pa_i(\sum\limits_{k=1}^3\bp_kU_{jk})=\sum\limits_{k=1}^3(\pa_i\bp_kU_{jk}+   \bp_k\pa_iU_{jk}),\quad \forall 1\leq i,j\leq 3. 
\end{equation*}
    \begin{equation*}
    \begin{split}
        \beta (Q(\bx))&= \sqrt{6}(\lam_1^3+\lam_2^3+\lam_3^3)\\
    &=\sqrt{6}\left((\f{\sqrt{6}}{6}+s)^3+ (\f{\sqrt{6}}{6}-s-\f{\sqrt{6}}{3}s^2+O(s^3))^3+ (-\frac{\sqrt{6}}{3}+\f{\sqrt{6}}{3}s^2+O(s^3))^3\right)\\
        &=-1+9s^2+\zeta(s)\\
        &=-1+\f{9}{2}|U|^2+\zeta(s),
    \end{split}
    \end{equation*}
for some function $\zeta(s)\sim O(s^3)$.

Consequently, on a small neighborhood $B_r(\bx_0)$, we can define the following energy functional
for $U:B_r(\bx_0)\to\mathcal{U}_{\bp}$:
\begin{equation}\label{ene: U}
    \begin{split}
        \mathcal{E}(U,B_r(\bx_0)):&=\int_{B_r(\bx_0)}\bigg\{\f12|\na 
 U|^2+\sqrt{6}\sum\limits_{1\leq i,j,k\leq 3} \pa_i \bp_j\,\pa_i \bp_k \,U_{jk}\\
  &-\f{\sqrt{6}}{4}|U|^2+\zeta(s)+(\na Q-\na R):\na R+\f12|\na R|^2\bigg\}\,d\bx,
    \end{split}
\end{equation}
Here $\na R$, $\na Q-\na R$ can be expressed in terms of $U$ and $\bp$. More specifically, the last two terms in \eqref{ene: U} can be written as follows: 
\begin{equation}\label{ene:higher order terms}
    \begin{split}
        &\int_{B_r(\bx_0)} \left\{ (\na Q-\na R):\na R+\f12|\na R|^2 
 \right\}\,d\bx\\
 =&\int_{B_r(\bx_0)} \Bigg\{ \bigg[ -\f{\sqrt{6}}{2}(\na \bp\otimes \bp+ \bp\otimes \na \bp)+\na U\bigg]: \bigg[\left(\f{\sqrt{6}}{3}+\f{d}{ds}\left(\f{\tau(s)}{s}\right)\right)\na s\left(\f{U}{2}+\f{3s}{2}(\bp(\bx)\otimes \bp(\bx)-\f13\mathrm{Id})\right)\\
 &\qquad\qquad  + \left(\f{\sqrt{6}}{3}s+\f{\tau(s)}{s}\right) \left(  \f{\na U}{2}+\f32 \na s(\bp\otimes\bp-\f13\mathrm{Id})+\f{3s}{2}(\na\bp\otimes \bp+\bp\otimes\na\bp) \right)\bigg]\\
 &\qquad + \left(\f{\sqrt{6}}{3}s+\f{\tau(s)}{s}\right)^2\left(\f{|\na U|^2}{4}+\f32|\na s|^2+\f{9s^2}{2}|\na\bp|^2-3s\sum\limits_{i,j,k}\pa_i\bp_j\pa_i\bp_k U_{jk}\right)\\
 &\qquad + 2\left(\f{\sqrt{6}}{3}+\f{d}{ds}\left(\f{\tau(s)}{s}\right)\right)^2s^2|\na s|^2\Bigg\}\,d\bx\\
 =&\int_{B_r(\bx_0)} \bigg\{  \sqrt{6} s
 |\na s|^2 -3s^2 |\na \bp|^2+s\left(\sum_{i,j,k}\pa_i\bp_j\pa_i\bp_kU_{jk}\right)+\f{\sqrt{6}}{6}s|\na U|^2 \bigg\}  \,d\bx+\mathcal{R}(B_r,\bp,U),
    \end{split}
\end{equation}
where 
\begin{equation}\label{est:higher order terms R}
    \mathcal{R}(B_r(\bx_0),\bp,U):= \int_{B_r(\bx_0)} \left(c_1(s)|\na s|^2+ c_2(s)|\na \bp|^2+c_3(s)\sum_{i,j,k}\pa_i\bp_j\pa_i\bp_kU_{jk}+c_4(s)|\na U|^2\right)\,d\bx,
\end{equation}
for quantities $c_1(s),c_3(s), c_4(s)\sim O(s^2)$, $c_2(s)\sim O(s^3)$.

\eqref{ene: U} together with \eqref{ene:higher order terms} and \eqref{est:higher order terms R} 
imply that

    \begin{equation}\label{ene: s,n leading term}
    \begin{split}
        \mathcal{E}(U,B_r(\bx_0))=&\int_{B_r} \bigg\{\left(\f12+\f{\sqrt{6}}{6}s+c_4(s)\right)|\na 
 U|^2+\left(\sqrt{6}+s+c_3(s)\right)\sum\limits_{1\leq i,j,k\leq 3} \pa_i \bp_j\,\pa_i \bp_k \,U_{jk}\\
  &-\f{\sqrt{6}}{4}|U|^2+\zeta(s)+\left(\sqrt{6}s+c_1(s)\right)|\na s|^2+\left(-3s^2+c_2(s)\right)|\na \bp|^2
 \bigg\}\,d\bx.\\
    \end{split}
\end{equation}
{Since $s=\f{|U|}{2}$ and $c_i(s)$ are explicit in $s$, they can all be regarded as functions of $U$}. Then it is not hard to see that $U\in H^1(B_r(\bx_0),\mathcal{U}_{\bp})$ is a minimizer of
$\mathcal{E}(\cdot, B_r(\bx_0))$ with respect to the Dirichlet boundary condition on $\pa B_r(\bx_0)$:
\begin{equation}\label{minimality}
\mathcal{E}(U, B_r(\bx_0))\le\mathcal{E}(V, B_r(\bx_0)),
\ \forall V\in H^1(B_r(\bx_0), \mathcal{U}_{\bp}),
\ {\rm{with}}\ V=U \ {\rm{on}}\ \pa B_r(\bx_0).
\end{equation}

For small $r$, we denote $k=k(\bx_0)$ and define the following rescaled functions
\begin{equation}\label{rescaling}
    \bp_r(\bx):=\bp(\bx_0+r\bx),\ \ U_r(\bx):=\f{1}{r^{k}}U(\bx_0+r\bx),\ \ s_r(\bx):=\f{1}{r^{k}}s(\bx_0+r\bx),\quad \bx\in B_1.
\end{equation}
By \eqref{polynomial repres}, 
\begin{equation}\label{Ur converge}
    \lim\limits_{r\ri 0}U_r(\bx)=\bar{U}(\bx), \quad \bx\in B_1,
\end{equation}
where $\bar{U}(\bx)$ is a matrix-valued function, whose components $\bar{U}_{ij}$ are all degree $k$ homogeneous polynomials. It follows from \eqref{minimality} that $U_r(\bx)$ minimizes the following energy functional:
\begin{equation}\label{ene:Ur}
\begin{split}
    &\mathcal{E}_r(U_r)\\
    &:=\int_{B_1}\bigg\{ \left(\f12+\f{\sqrt{6}}{6}r^{k}s+r^{2k}\hat{c}_4(s_r)\right)|\na 
 U_r|^2 + \f{\sqrt{6}+r^{k}s_r+r^{2k}\hat{c}_3(s_r)}{r^{k}} \sum\limits_{1\leq i,j,l\leq 3} \pa_i \bp_{r,j}\,\pa_i \bp_{r,l} \,U_{r,jl}\\
 &\qquad\quad-\f{\sqrt{6}r^2}{4}|U_r|^2+r^{k+2}\zeta(s_r)+r^{k}(s_r+r^{k}\hat{c}_1(s_r))|\na s_r|^2\\
 &\qquad\quad +(-3s_r^2+r^{k}\hat{c}_2(s_r))|\na \bp_r|^2\bigg\}\,d\bx\\
 &=: \mathcal{E}_1(U_r)+ \mathcal{E}_2(U_r)+ \mathcal{E}_3(U_r),
 \end{split}
\end{equation}
{where $\hat{c}_i(s_r(\bx)):=\f{c_i(s(\bx_0+r\bx))}{r^{2k}}\sim O(1)$ for $i=1,3,4$ and $\hat{c}_2(s_r(\bx)):=\f{c_2(s(\bx_0+r\bx))}{r^{3k}}\sim O(1)$.} And $\mathcal{E}_1(U_r)$, $\mathcal{E}_2(U_r)$ and $\mathcal{E}_3(U_r)$ are defined as follows: 
\begin{align}
\label{def:E1}   \mathcal{E}_1(U_r):= & \int_{B_1} \f12|\na U_r|^2\,\,d\bx,\\ 
\label{def:E2}     \mathcal{E}_2(U_r):= & \int_{B_1} \f{\sqrt{6}}{r^{k}}  \sum\limits_{1\leq i,j,l\leq 3} \pa_i \bp_{r,j}\,\pa_i \bp_{r,l} \,U_{r,jl}\,d\bx,\\
    \label{def:E3}\mathcal{E}_3(U_r):= & \int_{B_1} \bigg\{\left(\f{\sqrt{6}r^ks_r}{6}+r^{2k}\hat{c}_4(s_r)\right)|\na U_r|^2\\
    \nonumber &\qquad +\left(s_r+r^{k}\hat{c}_3(s_r)\right)  \sum\limits_{1\leq i,j,l\leq 3} \pa_i \bp_{r,j}\,\pa_i \bp_{r,l} \,U_{r,jl}-\f{\sqrt{6}r^2}{4}|U_r|^2\\
   \nonumber  &\qquad +r^{k+2}\zeta(s_r)+r^{k}(s_r+r^{k}\hat{c}_1(s_r))|\na s_r|^2 +(-3s_r^2+r^{k}\hat{c}_2(s_r))|\na \bp_r|^2\bigg\}\,d\bx.
\end{align}

The following lemma asserts that $\mathcal{E}_3$ is relatively small compared to $\mathcal{E}_1$.

\begin{lemma}\label{lemma: E3 small}
    There exists a constant $C=C(x_0)>0$ and $r_0>0$ such that 
    \begin{equation}\label{est: E3 small}
        \mathcal{E}_3(U_r)\leq Cr^{\min \{2,k\}}  \int_{B_1}  |\na U_r|^2\,d\bx,\quad \forall r\leq r_0. 
    \end{equation}
\end{lemma}

\begin{proof}
    Firstly, from  \eqref{Ur converge}, we have the bound 
    \begin{equation*}
        C_1\leq \int_{B_1}|\na U_r|^2\,d\bx\leq C_2, 
    \end{equation*}
    for some positive constants $C_1,C_2$. Therefore, it suffices to show that 
    \begin{equation}\label{est: bdd on E3}
        \mathcal{E}_3(U_r)\leq Cr^{\min\{2,k\}}, \quad \text{for some }C>0.
    \end{equation}
By the definitions of $s_r$ and $U_r$, we have $s_r(\bx)\sim O(1)$ and $U_r(x)\sim O(1)$. It remains to estimate the terms involving $\na \bp_r$. Indeed, since $\bp$ is analytic, we have 
\begin{equation*}
    |\na \bp_r(\bx)|=r|D\bp(\bx_0+r\bx)|\leq r \|\bp\|_{C^1(\bar{B}_1)} .
\end{equation*}
Then \eqref{est: bdd on E3} follows directly by estimating each term individually. 
\end{proof}
\section{Classification of local behaviors of $U(\bx)$} \label{sec:classification}

In this section we will present our classification of local behavior of $U(\bx)$ near the defect set $\mathcal{S}(Q)$  through Proposition \ref{prop: k=1 profile}, \ref{prop: rectifiable S1} and \ref{prop: U, k>=2}. We first establish an upper bound for the vanishing order $k(\bx)$ of $U(\bx)$ for $\bx\in \mathcal{S}(Q)$ using a compactness argument.

\begin{lemma}\label{lemma:upper bdd k}
    There exists a positive integer $N$, depending only on $\|Q\|_{C^1(\bar{B}_1)}$ and the constant $\delta$ in \eqref{bdy cond},  such that for any $\bx\in \mathcal{S}(Q)$, 
    \begin{equation*}
    k(\bx)\leq N.
    \end{equation*}
\end{lemma}
\begin{proof}
We prove it by contradiction. Suppose there were a sequence of minimizing maps $\{Q_n\}_{n=1}^\infty\in C^{\omega}(\bar{B}_1,\mathbb{S}^4)$ of $\mathcal{E}(Q,B_1)$ such that 
\begin{align}
  \label{Qn: uniform H1 bdd}  &\quad \|Q_n\|_{C^1(\bar{B}_1,\mathbb{S}^4)} \leq M,\quad \forall n\in\mathbb{N}^+,\\
   \label{Qn: bdy cond} &\beta(Q_n(\bx)) >-1+\delta,\quad \text{for some }\delta>0,\ \forall \bx\in\pa B_1, \ n\in\mathbb{N}^+,
\end{align}
and a sequence of $\bx_n\in B_1\cap \mathcal{S}(Q_n)$ such that 
\begin{equation*}
   k(\bx_n)=n,\quad \forall n\in\mathbb{N}^+. 
\end{equation*}
From the boundary condition \eqref{Qn: bdy cond} and the uniform bound on $|\na Q_n(\bx)|$, we obtain
\begin{equation*}
    \dist(\bx_n,\pa B_1)>\eta,\quad \text{for some }\eta=\eta(\delta,M)>0,\quad \forall n\in\mathbb{N}^+.
\end{equation*}

By the compactness of minimizing maps of $\mathcal{E}(Q,B_1)$ in $H^1(B_1,\mathbb{S}^4)$, we conclude that, up to a subsequence which is still denoted by $\{Q_n\}_{n=1}^\infty$,
\begin{align*}
    &Q_n\ri \tilde{Q}\  \text{ in } C^l_{loc}(B_1,\mathbb{S}^4),\  \forall l\in \ \mathbb{N}^+,\ \text{ uniformly on }\bar{B_1};\\
    & \bx_n\ri \tilde{\bx} \quad \text{for some }\tilde{\bx}\in B_1 \text{ such that } \dist(\tilde{\bx},\pa B_1)\geq \eta,
\end{align*}
where $\tilde{Q}\in C^{\omega}(B_1,\mathbb{S}^4)\cap W^{1,\infty}(\bar{B}_1)$ is a minimizing map
of $\mathcal{E}(Q, B_1)$ with $\tilde{\bx}\in\mathcal{S}(\tilde{Q})\cap B_{1-\f{\eta}{2}}$. Since $\tilde{U}$ has infinite vanishing order at $\tilde{\bx}$,  the unique continuation property implies that the analytic function $\tilde{Q}+\f{\sqrt{6}}{2}(\tilde{\bp}\otimes \tilde{\bp}-\f13\mathrm{Id})$ is identically zero. Consequently,  $\tilde{Q}$ must be uniaxial everywhere and $\beta(\tilde{Q}(\bx))\equiv -1$ in $B_1$. This yields a contradiction with $\beta(\tilde{Q}(\bx))\geq -1+\delta$ for $\bx\in \pa B_1$.
\end{proof}

Lemma \ref{lemma:upper bdd k} implies the following decomposition of $\mathcal{S}(Q)$ in terms of the vanishing order $k$:
\begin{equation*}
    \mathcal{S}(Q)=\displaystyle\bigcup_{j=1}^N \mathcal{S}_j(Q),\quad \text{where } \mathcal{S}_j(Q):=
    \big\{\bx\in\mathcal{S}(Q):k(\bx)=j\big\}.
\end{equation*}

We start with the structure of $\mathcal{S}_1(Q)$ and local profile of $U$ near $\mathcal{S}_1(Q)$.  Take $\bx_0\in \mathcal{S}_1(Q)$. Recall that the eigenframe $\{\bn,\bm,\bp\}$ and the leading eigenvector $\bn$ are discontinuous at $\bx_0$.   Let $\bar{U}(\cdot)$ denote the leading polynomial at $\bx_0$, i.e.
\begin{equation*}
    U(\bx_0+\bx)=\bar{U}(\bx)+O(|\bx|^2),\quad \text{ for }|\bx|\ll 1.
\end{equation*}
Since the vanishing order $k(\bx_0)=1$, each non-zero component of $\bar{U}$ is a linear function.

\begin{prop}\label{prop: k=1 profile}
    Suppose $\bx_0\in \mathcal{S}_1(Q)$ and $\bar{U}(\bx)$ is the blow-up limit of $U$ at $\bx_0$. There exists a unit vector $\mathbf{e}\in\mathbb{S}^2$ such that  the following properties hold:
\begin{enumerate}
\item The nodal set of $\bar{U}$, i.e. $\{\bx:\bar{U}(\bx)=0\}$, is the straight line $l_{\mathbf{e}}:=\{\bx: \bx=c\,\mathbf{e},\,c\in\mathbb{R}\}$.  
\item$\bar{U}$ is invariant along the direction $\mathbf{e}$, and therefore can be considered as a function of two variables. 
\item The winding number of the leading eigenvector $\bn$ of $\bar{U}$ around $l_{\mathbf{e}}$ is $\pm\f12$.
\end{enumerate}
\end{prop}

\begin{proof}
    Upon a possible rotation of the coordinate, we assume $\bp(\bx_0)=(0,0,1)$. Note that for any $\bx\in B_1$, $\bp_r(\bx)$ is an eigenvector of $U_r(\bx)$ associated with the eigenvalue $0$.  As $r$ tends to 0, $U_r$ and $\bp_r$ converges to $\bar{U}(\bx)$, $\bp(\bx_0)$ uniformly, respectively. Consequently, it holds
    \begin{equation*}
        \bar{U}(\bx) \bp(\bx_0)=0,\quad \forall \bx\in B_1. 
    \end{equation*}
    Since $\bar{U}$ is traceless, symmetric and linear, there exist two vectors $\mathbf{a},\ \mathbf{b}\in\mathbb{R}^3$ such that 
    \begin{equation}\label{blowuplimit}
        \bar{U}(\bx)=\begin{bmatrix}
            \mathbf{a}\cdot \bx & \mathbf{b}\cdot \bx&0\\
            \mathbf{b}\cdot \bx & -\mathbf{a}\cdot \bx&0\\
            0&0&0
        \end{bmatrix}.
    \end{equation}

\begin{figure}[htb!]
    \centering
    \begin{subfigure}[t]{0.5\textwidth}
        \centering


\tdplotsetmaincoords{70}{120}

\begin{tikzpicture}[tdplot_main_coords, scale=1.5]

\draw[->, thick] (0,0,0) -- (1.8,0,0) ;
\draw[->, thick] (0,0,0) -- (0,1.8,0) ;
\draw[->, thick] (0,0,-1.5) -- (0,0,2.4) node[anchor=north east] {$\mathbf{a}$};
\draw[->, thick] (0,1.5,1.8) -- ++ (0.6,0,0)node[anchor=north east] {$\mathbf{n}_0(\lambda)$} ;
\draw[->, thick] (0,1.5,1.8) -- ++ (0,0.6,0)
node[anchor=north] {$\mathbf{n}_0(\lambda)\times \bp$};
\draw[->, thick] (0,1.5,1.8) -- ++ (0,0,0.5) node[anchor=north east] {$\mathbf{p}$};

\foreach \z in {-1.7,-0.7,0.7,1.7} {
  \foreach \x in {-1.2,-0.6,0,0.6,1.2} {
    \foreach \y in {-1.2,-0.6,0,0.6,1.2} {
      \pgfmathsetmacro{\u}{\z}
      \pgfmathsetmacro{\uz}{0.23*\u}
      \ifdim\z pt > 0pt
        \ifdim\z pt=0.7pt
        \draw[-, thick, blue] (\x,\y,\z) -- ++(\uz,0,0);
        \else
        \draw[-, thick, blue!50, opacity=0.6] (\x,\y,\z) -- ++(\uz,0,0);
        \fi
      \else
      \ifdim\z pt=-0.7pt
      \draw[-, thick, red] (\x,\y,\z) -- ++(0,-\uz,0);
      \else
       \draw[-, thick, red!50, opacity=0.6] (\x,\y,\z) -- ++(0,-\uz,0);
       \fi
    
      \fi
    }
  }
}

\end{tikzpicture}

        \captionsetup{width=0.8\textwidth}
        \caption{Case 1:  the eigenframe remains unchanged across $\mathbf{a}\cdot \mathbf{x}=0$, with only the leading eigenvector changing from $\mathbf{n}_0$ to $\mathbf{n}_0\times \mathbf{p}$.}
        \label{Fig:case1}
    \end{subfigure}%
    ~ 
    \begin{subfigure}[t]{0.5\textwidth}
        \centering


\tdplotsetmaincoords{70}{120}
\begin{tikzpicture}[tdplot_main_coords, scale=1.8]

\draw[->, thick] (0,0,0) -- (1.8,0,0) node[anchor=north east] {$\mathbf{a}$};
\draw[->, thick] (0,0,0) -- (0,1.8,0) node[anchor=north west] {$\mathbf{b}$};
\draw[->, thick] (0,0,-1) -- (0,0,1.8) node[anchor=south] {$\mathbf{a} \times \mathbf{b}$};

\foreach \z in {-1,0,1} {
  \foreach \x in {-1.2,-0.7,-0.2,0.2,0.7,1.2} {
    \foreach \y in {-1.2,-0.7,-0.2,0.2,0.7,1.2} {
      \pgfmathsetmacro{\u}{\x}
      \pgfmathsetmacro{\v}{\y}
      \pgfmathsetmacro{\norm}{sqrt(\u*\u + \v*\v)}
      \pgfmathsetmacro{\ux}{0.35*(sqrt((\norm+\x)/2))}
      \pgfmathsetmacro{\vy}{0.35*(sqrt((\norm-\x)/2)* \y/ sqrt(\y*\y))}
      \ifdim\z pt=0pt
        \draw[-, thick, blue] (\x,\y,0) -- ++(\ux,\vy,0);
      \else
        \draw[-, thick, blue!40, opacity=0.5] (\x,\y,\z) -- ++(\ux,\vy,0);
      \fi
    }
  }
}

\end{tikzpicture}

    
        \captionsetup{width=0.8\textwidth}
        \caption{Case 2: the leading eigenvector $\bn(x)$ forms a $\f12$ degree defect around $\mathbf{a}\times\mathbf{b}$.}
        \label{Fig:case2}
    \end{subfigure}
    \caption{Vector field $s(\bx)\bn(\bx)$ for $k(\bx_0)=1$}
    \label{Fig:k=1}
\end{figure}

We divide the discussion into two cases:
\begin{case}
$\mathbf{a}\times \mathbf{b}=\mathbf{0}$. We assume $|\mathbf{a}|\neq 0$ and $\mathbf{b}=\lambda \mathbf{a}$ for some $\lambda\in\mathbb{R}$.  Then the nodal set  $\{|\bar{U}|=0\}$ can be identified with a two-dimensional plane $\{\bx: \mathbf{a}\cdot \bx=0\}$. More specifically, the positive eigenvalue $s(\bx)$ and its corresponding eigenvector $\bn(\bx)$ are given by
\begin{equation*}
    s(\bx)=|\mathbf{a}\cdot\bx|\sqrt{1+\lambda^2},\quad\bn(\bx)=\begin{cases}(\f{\lambda}{\sqrt{2+2\lambda^2-2\sqrt{1+\lambda^2}}}, \f{\sqrt{1+\lam^2}-1}{\lam\sqrt{2+2\lambda^2-2\sqrt{1+\lambda^2}}},0), &  \mathbf{a}\cdot\bx>0,\\
    (\f{\sqrt{1+\lam^2}-1}{\lam\sqrt{2+2\lambda^2-2\sqrt{1+\lambda^2}}},-\f{\lambda}{\sqrt{2+2\lambda^2-2\sqrt{1+\lambda^2}}},0), &  \mathbf{a}\cdot\bx<0.
    \end{cases}
\end{equation*}
This implies that $\bn$ equals to a constant vector $\bn_0(\lambda)$ on $\{\mathbf{a}\cdot\bx>0\}$, while on $\{\mathbf{a}\cdot\bx<0\}$, it becomes $\bn_0(\lambda)\times \bp(\bx_0)$.  In this case, $\{|\bar{U}|=0\}$ is where the exchange of eigenvalues occurs. As $\bx$ passes through the plane, the leading eigenvector shifts from $\bn_0$ to its perpendicular vector $\bn_0\times \bp(x_0)$.  Therefore, $\bx_0$ cannot be considered as a defect point by definition \eqref{def:defect}, as the eigenframe remains continuous in its neighborhood. 
\end{case}
\begin{case}
    $\mathbf{a}\times \mathbf{b}\neq \mathbf{0}$.  Then $\{\bar{U}=0\}$ forms a straight line $l=\{\bx=c\,\mathbf{a}\times \mathbf{b},\,c\in\mathbb{R}\}$, and $\bar{U}$ is invariant along the direction $\mathbf{a}\times \mathbf{b}$. Therefore, $\bar{U}$ can be essentially considered as a function of two variables, corresponding to coordinates of two perpendicular directions of the axis $l$. For $\bx(\theta):=\cos\theta \f{\mathbf{a}}{|\mathbf{a}|^2}+\sin\theta\f{\mathbf{b}}{|\mathbf{b}|^2}$, when $\theta$ varies from $0$ to $2\pi$, $\bx(\theta)$ traces a full circle around the axis $l$ .  The leading eigenvector $\bn(\bx)$ is $(\cos{\f{\theta}{2}},\sin{\f{\theta}{2}},0)$, indicating that the winding number of $\bn$ around $l$ is $\pm\f12$. 

\end{case}
    
 The proof is complete.
\end{proof}

One can apply the dimension reduction principle as carried out in \cite{lin1991nodal} to derive that the Hausdorff dimension of the set $\mathcal{S}_1(Q)$ is $1$. Indeed, we can show that it is a countably $1$-rectifiable set.

\begin{prop}[Rectifiability] \label{prop: rectifiable S1}
 The set $\mathcal{S}_1(Q)$ is a countably $1$-rectifiable set.
\end{prop}

\begin{proof}
    Take any $\bx_0\in \mathcal{S}_1$. Let $l_{\mathbf{e}}=\{\bx=c\,\mathbf{e} \ \big|\ c\in\mathbb{R}\}$ be the invariant axis of the blow up limit $\bar{U}$ of $U$ at $\bx_0$. 
   We claim that  
    \begin{equation}\label{cone_angle}
        \lim\limits_{r\rightarrow 0}
        \sup_{\bx\in \mathcal{S}_1(Q)\cap B_r(\bx_0)}\frac{\mathrm{dist}(\bx,l_\mathbf{e})}{|\bx-\bx_0|}=0.
    \end{equation}
    To prove \eqref{cone_angle}, we argue by contradiction. If \eqref{cone_angle} were false, there would exist $c_0>0$, a radius sequence $r_i\to 0$, and a point sequence $\{\bx_i\}\subset\mathcal{S}_1(Q)\cap B_{r_i}(\bx_0)$ such that
\begin{equation}\label{cone_angle1}
 \frac{\mathrm{dist}(\bx_i,l_\mathbf{e})}{|\bx_i-\bx_0|}\ge c_0. 
 \end{equation}
 For $i\ge 1$, let $\bar{\bx}_i\in B_{r_i}(x_0)\cap l_{\bf{e}}$   be such that
 $$
 |\bx_i-\bar{\bx}_i|=\mathrm{dist}(\bx_i,l_\mathbf{e}).
 $$
 From \eqref{polynomial repres}, we have
 \begin{align*}
  0=U(\bx_i)=\bar{U}(\bx_i-\bx_0)
  +O(|\bx_i-\bx_0|^2).
 \end{align*}
Since $\bx_i-\bar{\bx}_i\perp l_{\bf e}$, it follows from \eqref{blowuplimit} that 
\[
\bar{U}(\bx_i-\bx_0)=\bar{U}(\bx_i-\bar{\bx}_i)
\]
so that
$$
\big|\bar{U}(\bx_i-\bx_0)\big|
\ge \sqrt{\bf{a}^2+\bf{b}^2} \big|\bx_i-\bar{\bx}_i\big|.
$$ 
Thus, we arrive at
\[
\mathrm{dist}(\bx_i,l_\mathbf{e})=\big|\bx_i-\bar{\bx}_i\big|=O(|\bx_i-\bx_0|^2).
\]
This contradicts \eqref{cone_angle1}.
Hence \eqref{cone_angle} holds.
It is well known that \eqref{cone_angle}
implies that for every point $\bx_0\in \mathcal{S}_1(Q)$, there exists a unique $1$-dimensional tangent line of $\mathcal{S}_1(Q)$ at $\bx_0$. 
Thus, by the classic geometric measure theory, $\mathcal{S}_1(Q)$ is a countably $1$-rectifiable set.
\end{proof}

\begin{rmk}
    Proposition \ref{prop: k=1 profile} and Proposition \ref{prop: rectifiable S1} demonstrate the half-degree profile of the leading eigenvector of $Q$ around disclination lines under the Lyuksyutov constraint. The asymptotic profiles of eigenvalues and their corresponding eigenvectors are similar to the characterizations of biaxial-ring defects in \cite{yu2020,ty}. However, we point out that in our analysis, the invariant axis for $\bar{U}$ does not coincide with the eigenvector $\bp(\bx_0)$. In contrast,  at the biaxial-ring defect core in \cite{yu2020,ty}, the tangent direction of the ring defect serves as the eigenvector associated with the negative eigenvalue,  due to the imposed symmetry hypotheses. 
\end{rmk}

For points with higher vanishing order $k\geq 2$, we begin by deriving the energy functional that the local blowup $\bar{U}(\bx)$ is expected to minimize, based on the computations carried out in Section \ref{sec:functional}. Lemma \ref{lemma: E3 small} asserts that $\mathcal{E}_3(U_r)$ is negligible compared to $\mathcal{E}_{1}(U_r)$. Therefore, to identify the leading-order term of the energy as $r\rightarrow 0$, it suffices to compare $\mathcal{E}_2$ with $\mathcal{E}_1$. 

Since $\bp$ is real analytic, we can write the Taylor expansion of $\na \bp_r$ up to degree $k-1$:
\begin{align*}
    \na \bp_r(\bx)&=r Dp(\bx_0+r\bx)\\
    &=r\Big(\sum\limits_{|\alpha|\leq k-2} \f{D\circ D^\alpha \bp(\bx_0)}{\alpha!}\bx^\alpha r^{|\al|}+r^{k-1} P_k(\bx) \Big)\\
    &=\sum\limits_{|\alpha|\leq k-2} \f{D\circ D^\alpha \bp(\bx_0)}{\alpha!}\bx^\alpha r^{|\al|+1}+r^{k} P_k(\bx),
\end{align*}
where $\alpha=(\alpha_1,\alpha_2,\alpha_3)$ is the multi-index with $|\alpha|=|\alpha_1|+|\alpha_2|+|\alpha_3|$ and $\alpha!=\alpha_1!\alpha_2!\alpha_3!$.  $P_k(\bx)$ is the reminder term, which is approximately a $(k-1)$-degree polynomial of $\bx$. 

It follows that  
\begin{align*}
    &\sum\limits_{1\leq i,j,l\leq 3} \pa_i \bp_{r,j}\,\pa_i \bp_{r,l} \,U_{r,jl}\\
    =&\sum\limits_{ i,j,l} \Big( \sum\limits_{|\alpha|\leq k-2} \f{D_i D^\alpha \bp_j(x_0)}{\alpha!}\bx^\alpha r^{|\al|+1} \Big) \cdot\Big(\sum\limits_{|\beta|\leq k-2} \f{D_i D^\beta \bp_l(x_0)}{\beta!}\bx^\beta r^{|\beta|+1}\Big) \cdot U_{r,jl}(\bx)+O(r^{k+1})\\
    =&\sum\limits_{i,j,l}\sum\limits_{m=2}^k \sum\limits_{|\al|+|\beta|=m-2} \Big(\f{D_i D^\al \bp_j(x_0)}{\al!}\cdot \f{D_i D^\beta \bp_l(x_0)}{\beta!} \Big)\bx^{\al+\beta} r^m \, U_{r,jl}(\bx) +O(r^{k+1})
\end{align*}

Let $V_m(\bx)$ denote the tensor-valued function, 
\begin{equation}\label{def: Vm}
    V_m(\bx):=\sum\limits_{i=1}^3\sum\limits_{|\al|+|\beta|=m-2} \f{D_i D^\al \bp(x_0)}{\al!}\otimes \f{D_i D^\beta \bp(x_0)}{\beta!} \bx^{\al+\beta},\quad 2\le m\le k, 
\end{equation}
where $V_m$ is determined by $D\bp(\bx_0)$, $D^2\bp(\bx_0)$,..., $D^{m-1}\bp(\bx_0)$. Substituting into the formula \eqref{ene:Ur} of the energy functional gives
\begin{equation}\label{formula: Er using Vm}
    \qquad \mathcal{E}_r(U_r)=\int_{B_1}\f12|\na U_r|^2\,d\bx +\f{\sqrt{6}}{r^k} \sum\limits_{m=2}^k r^{m} \int_{B_1} (V_m:U_r)\,d\bx+O(r)+\mathcal{E}_3(U_r).
\end{equation}

\begin{lemma}\label{lemma: Vm=0 m<k}
    For any $2\leq m<k$ and $\bx\in B_1$, $V_m(\bx)$ satisfies
    \begin{equation}\label{Vm:Q=0}
        V_m(\bx): Q=0, \quad \forall\,Q\in \mathcal{U}_{\bp(\bx_0)},
    \end{equation}
where $\mathcal{U}_{\bp(\bx_0)}:=\{U\in \mathcal{Q}_0:\, U\bp(\bx_0)=\mathbf{0}\}$.
\end{lemma}

\begin{proof}
Suppose \eqref{Vm:Q=0} were false for some $2\leq m<k$. We take $m$ as the smallest integer such that  \eqref{Vm:Q=0} does not hold. Denote $\bp_0=\bp(x_0)$ and define 
    \begin{equation}\label{def: Ym}
    \begin{split}
        Y_m(\bx):= &V_m-V_m(\bp_0\otimes\bp_0)-(\bp_0\otimes\bp_0)V_m+(V_m:\bp_0\otimes\bp_0)\bp_0\otimes\bp_0\\
        &\quad -\f{V_m:(\mathrm{Id}-\bp_0\otimes\bp_0)}{2}(\mathrm{Id}-\bp_0\otimes\bp_0).
        \end{split}
    \end{equation}
{$Y_m(\bx)$ is the orthogonal projection of $V_m(\bx)$ onto the space $\mathcal{U}_{\bp_0}$, the term $-\f{V_m:(\mathrm{Id}-\bp_0\otimes\bp_0)}{2}(\mathrm{Id}-\bp_0\otimes\bp_0)$  is the Lagrange multiplier with the traceless constraint.} It follows that \eqref{Vm:Q=0} is equivalent to $Y_m(\bx)=0$ for all $\bx\in B_1$.
  
    Take $c>0$ and $\rho\in(0,1)$ as constants to be determined later. We define the following energy competitor:
    \begin{equation*}
        Z_r(\bx):=\begin{cases}
            -cY_m(\bx), & |x|\leq \rho,\\
            -cY_m(\bx)\f{1-|x|}{1-\rho}+U_r(\bx)\f{|x|-\rho}{1-\rho}, & \rho<|x|\leq 1.
        \end{cases}
    \end{equation*}
Direct calculation implies
\begin{align*}
    \int_{B_1} V_m(\bx):Z_r(\bx)\,d\bx&=-c\int_{B_\rho} V_m:Y_m\,d\bx+\int_{B_1\setminus B_\rho} (  -cY_m(\bx)\f{1-|x|}{1-\rho}+U_r(\bx)\f{|x|-\rho}{1-\rho}):V_m \,d\bx\\
    &=\leq -c\int_{B_\rho} |Y_m|^2\,d\bx +c\int_{B_1\setminus B_\rho}|Y_m|^2\,d\bx+\int_{B_1\setminus B_\rho}|U_r:V_m|\,d\bx.
\end{align*}
Since $\int_{B_1}|Y_m|^2\,d\bx>0$, we can always find suitable $c,\rho$ such that
\begin{equation*}
    \int_{B_1} V_m:Z_r\,d\bx\leq -1+\int_{B_1} V_m:U_r\,d\bx,  
\end{equation*}
for any $r$ sufficiently small.  Fix this choice of $c$ and $\rho$, from \eqref{formula: Er using Vm} we obtain
\begin{equation*}
    \mathcal{E}_r(U_r)-\mathcal{E}_r(Z_r) \geq  \f{\sqrt{6}}{r^{k-m}} + \int_{B_1} \f12(|\na U_r|^2-|\na Z_r|^2)\,d\bx +o(\f{1}{r^{k-m}})>\f{\sqrt{6}}{2r^{k-m}},
\end{equation*}
for any $r$ sufficiently small, here we used the fact that $U_r$ and $Z_r$ are uniformly bounded in $H^1(B_1)$ regardless of $r$.  This contradicts to the fact that $U_r$ minimizes $\mathcal{E}_r$.
\end{proof}

It follows from Lemma \ref{lemma: Vm=0 m<k} and Lemma \ref{lemma: E3 small} that the functional \eqref{formula: Er using Vm} reduces to 
\begin{equation}\label{formula: Er, only Vk}
    \mathcal{E}_r(U_r)=\int_{B_1} \f12|\na U_r|^2\,d\bx+\sqrt{6}\int_{B_1} (V_k:U_r) \,d\bx +O(r),
\end{equation}
 where $V_k(\bx)$ is the polynomial defined in \eqref{def: Vm}, with $m$ replaced by $k$.  In the limiting $r\ri 0$ case, we have the following characterization on the structure of $\bar{U}$ and $\mathcal{S}_k(Q)$, which completes our classification of the set $\mathcal{S}(Q)$.

 \begin{prop}\label{prop: U, k>=2}
     For any $k\geq 2$, $\mathcal{S}_k(Q)$ consists of isolated points and a countable union of $C^1$ curves. Suppose $\bx_0\in\mathcal{S}_k(Q)$ and let $\bar{U}(\bx)$ denote the blow-up limit of  $U$ at $\bx_0$.  Then $\bar{U}$ minimizes the energy
     \begin{equation}\label{def: Ek}
         \mathcal{E}^k(U):= \int_{B_1}\Big( \f12|\na  U|^2 +\sqrt{6} V_k: U\Big)\,d\bx, 
     \end{equation}
   in $H^1(B_1, \mathcal{U}_{\bp(\bx_0)})$
    subject to the Dirichlet boundary condition $U=\bar{U}$ on $\partial B_1$.  And $\bar{U}$ solves the system
     \begin{equation}\label{EL eq for k>=2}
         \Delta \bar{U}= \sqrt{6}Y_k,
     \end{equation}
     where $Y_k$ is defined in \eqref{def: Ym}, with $m$ replaced by $k$. There exists a unit vector $\mathbf{e}\in\mathbb{S}^2$ such that the nodal set of $\bar{U}$ is the line $l_{\mathbf{e}}=\{\bx=c\mathbf{e}\ \big|\ c\in\mathbb{R}\}$. Moreover, $\bar{U}$ can be viewed as a function of two variables, with $\mathbf{e}$ as the invariant direction. 
 \end{prop}
\begin{proof}
    We first prove the minimizing property of $\bar{U}$. Let $X\in H_0^1(B_1,\mathcal{U}_{\bp(\bx_0)})$. Since $\bp_r(\bx)\ri \bp(\bx_0)$ in $C^1(B_1)$, it is clear that one can find a function $X_r\in H_0^1(B_1,\mathcal{U}_{\bp_r})$ for any $r$,  where $\mathcal{U}_{\bp_r}$ is defined as in \eqref{target U},  with $\bp$ replaced by $\bp_r$, such that $X_r\to X$ strongly in $H^1(B_1)$. The minimality of $U_r$ implies
    \begin{align*}
        0&\leq \mathcal{E}_r(U_r+X_r)-\mathcal{E}_r(U_r)\\
       & = \int_{B_1} \left( \na U_r: \na X_r +\f12|\na X_r|^2 +\sqrt{6} V_k:X_r\right)\,d\bx+ O(r).
    \end{align*}
    Utilizing the strong $H^1$-convergence of $(X_r,U_r)$ to $(X,\bar{U})$, we obtain that, as $r\to 0$,
    \begin{equation*}
        0\leq\int_{B_1} \left( \na \bar{U}: \na X +\f12|\na X|^2 +\sqrt{6} V_k:X\right)\,d\bx=\mathcal{E}^k(\bar{U}+X)-\mathcal{E}^k(\bar{U}).
    \end{equation*}
    Consequently, $\bar{U}$ minimizes \eqref{def: Ek} and satisfies the Euler-Lagrange system \eqref{EL eq for k>=2}, as claimed. Furthermore, $\bar{U}$ is  the sum of a degree $k$ homogeneous polynomial determined by $V_k$ and a degree $k$ harmonic polynomial.

    Next, we establish the local symmetry of $\bar{U}$.  Set
\begin{equation*}
    \mathcal{S}_{k,\bx_0}:=\Big\{\bx\in B_1: \bar{U}(\bx)=0, \ \bar{U}(\bx+r\mathbf{y})=r^k\bar{U}(\bx+\mathbf{y}),\  \forall \bx+\mathbf{y}\in B_1, r\in (0,1)\Big\}.
\end{equation*}
as the subset of the nodal set of $\bar{U}$ that consists of degree $k$ homogeneous points. If $\mathcal{S}_{k,\bx_0}\backslash \{0\}\neq \varnothing$, we assume, up to a possible change of coordinates, that $(0,0,\f12)\in \mathcal{S}_{k,\bx_0}$. Since $\bar{U}$ is degree $k$ homogeneous with respect to both the origin and $(0,0,\f12)$, it is invariant along the $x_3$-direction. We claim that $\mathcal{S}_{k,\bx_0}=\{(0,0,x_3): x_3\in(-1,1)\}$. Suppose instead there is a point $\mathbf{z}$ off the $x_3$-axis that also belongs to $\mathcal{S}_{k,\bx_0}$. Then $\bar{U}$ would vanish on the plane spanned by $\mathbf{z}$ and $x_3$-axis, implying that each element of $\bar{U}$ is a degree $k$ polynomial depending on only one variable. Arguing as in Case 1 of the proof of Proposition \ref{prop: k=1 profile}, this yields a contradiction, as $\bx_0\not\in \mathcal{S}_{k,\bx_0}$ due to the continuity of the eigenframe around $\bx_0$.    

Let $\tilde{\mathcal{S}}_k$ denote the subset of $\mathcal{S}_{k,\bx_0}$, whose tangent map $\bar{U}$ depends on two variables. The above reasoning implies that $\tilde{\mathcal{S}}_k=\{\bx_0\in\mathcal{S}_k(Q):\, \mathcal{S}_{k,\bx_0}\setminus\{0\}\neq \varnothing\}$. The 1-rectifiability of $\tilde{\mathcal{S}}_k$ follows immediately from the same arguments in the proof of Proposition \ref{prop: rectifiable S1}.  We show that $\mathcal{S}_k\backslash \tilde{\mathcal{S}}_k$ only contains isolated points. Suppose, for contradiction, there exists a non-isolated point $\bx\in \mathcal{S}_k\backslash \tilde{\mathcal{S}}_k$. Then there exists a sequence of points $\{\bx_i\}\subset \mathcal{S}_k\backslash \tilde{\mathcal{S}}_k$ such that $\bx_i\ri \bx$. If we consider the rescaling \eqref{rescaling} around $\bx$ with $r_i=2|\bx_i-\bx|$, then $U_{r_i}$ converges in $C^l$ (for all $l\ge 1$) to a degree $k$ homogeneous polynomial $\bar{U}$ and $\f{\bx_i-\bx}{r_i}\ri \mathbf{y}$ for some point $\mathbf{y}\in B_1$ satisfying $|\mathbf{y}|=\f12$. The $C^l$ convergence of $\{U_{r_i}\}$ implies that $\mathbf{y}\in \mathcal{S}_{k,\bx}$ , thus rendering $\bar{U}$ invariant along the $\mathbf{y}$ direction. This contradicts with the assumption that $\bx\not\in \tilde{\mathcal{S}}_k$. The proof is complete.    

\end{proof}

\begin{rmk}
So far, all of our descriptions of the local structure near $\mathcal{S}(Q)$ apply only to minimizer $Q$ of \eqref{ene} subject to the Lyuksyutov constraint $|Q|=1$. However, this constraint can be relaxed, and the arguments can be generalized to minimizer
of the more general energy functional \eqref{def: E_{lambda,mu}} in the regime where $\mu\gg 1$. We assume the same boundary condition \eqref{bdy cond} and let $\mathcal{S}$ denote the defect set defined in \eqref{def:defect}. For $\bx_0\in \mathcal{S}$, we continue to define $U(\bx)$ as the leading term in the local Taylor expansion of $Q(\bx)+\f{\sqrt{6}}{2}(\bp(\bx)\otimes \bp(\bx)-\f13\mathrm{I_d})$ around $\bx_0$, and denote by $k(\bx_0)$ the vanishing order. Then the following argument holds:
for any $\mu \gg 1$, there exists $\varepsilon(\mu) $ such that 
\begin{enumerate}
    \item $\varepsilon(\mu)\to 0$ as $\mu\to \infty$
    \item For any $\bx_0\in \mathcal{S}$, there exists a $\bar{U}(\bx)$, whose structure is predicted by Proposition \ref{prop: k=1 profile} ($k(\bx_0)=1$) or Proposition \ref{prop: U, k>=2} ($k(\bx_0)\geq 2$), such that 
    \begin{equation*}
       \lim\limits_{r\ri 0} \sup\limits_{\bx\in B_1} \left|  \f{1}{r^{k(\bx_0)}}U(\bx_0+r\bx)-\bar{U}(\bx) \right|\leq \varepsilon(\mu). 
    \end{equation*}
\end{enumerate}
The proof relies only on a compactness argument and we omit the details. 
\end{rmk}

\bigskip
\noindent{\bf Acknowledgments:}  
The first author is partially supported by an AMS-Simons travel grant. 
The second author is partially supported by NSF grant 2101224 and a Simons travel grant.

\bigskip

\bibliographystyle{acm}
\bibliography{line}

\begin{thebibliography}{10}

\bibitem{alper2018rectifiability}
{\sc Alper, O.}
\newblock Rectifiability of line defects in liquid crystals with variable
  degree of orientation.
\newblock {\em Archive for Rational Mechanics and Analysis 228}, 1 (2018),
  309--339.

\bibitem{alper2017defects}
{\sc Alper, O., Hardt, R., and Lin, F.-H.}
\newblock Defects of liquid crystals with variable degree of orientation.
\newblock {\em Calculus of Variations and Partial Differential Equations 56\/}
  (2017), 1--32.

\bibitem{Ball}
{\sc Ball, J.~M.}
\newblock Mathematics and liquid crystals.
\newblock {\em Molecular Crystals and Liquid Crystals 647}, 1 (2017), 1--27.

\bibitem{bpp}
{\sc Bauman, P., Park, J., and Phillips, D.}
\newblock Analysis of nematic liquid crystals with disclination lines.
\newblock {\em Archive for Rational Mechanics and Analysis 205}, 3 (2012),
  795--826.

\bibitem{bbh93}
{\sc Bethuel, F., Brezis, H., and H{\'e}lein, F.}
\newblock Asymptotics for the minimization of a {Ginzburg-Landau} functional.
\newblock {\em Calculus of Variations and Partial Differential Equations 1}, 2
  (1993), 123--148.

\bibitem{bbh94}
{\sc Bethuel, F., Brezis, H., H{\'e}lein, F., et~al.}
\newblock {\em {Ginzburg-Landau} vortices}, vol.~13.
\newblock Springer, 1994.

\bibitem{bcl}
{\sc Brezis, H., Coron, J.-M., and Lieb, E.~H.}
\newblock Harmonic maps with defects.
\newblock {\em Communications in Mathematical Physics 107}, 4 (1986), 649--705.

\bibitem{bronshtein1979smoothness}
{\sc Bronshtein, M.~D.}
\newblock Smoothness of roots of polynomials depending on parameters.
\newblock {\em Siberian Mathematical Journal 20}, 3 (1979), 347--352.

\bibitem{canevari1}
{\sc Canevari, G.}
\newblock Biaxiality in the asymptotic analysis of a {2D} {Landau- de Gennes}
  model for liquid crystals.
\newblock {\em ESAIM: Control, Optimisation and Calculus of Variations 21}, 1
  (2015), 101--137.

\bibitem{canevari2}
{\sc Canevari, G.}
\newblock Line defects in the small elastic constant limit of a
  three-dimensional {Landau-de Gennes} model.
\newblock {\em Archive for Rational Mechanics and Analysis 223}, 2 (2017),
  591--676.

\bibitem{cl2017}
{\sc Contreras, A., and Lamy, X.}
\newblock Biaxial escape in nematics at low temperature.
\newblock {\em Journal of Functional Analysis 272}, 10 (2017), 3987--3997.

\bibitem{de1972types}
{\sc De~Gennes, P.}
\newblock Types of singularities permitted in the ordered phase.
\newblock {\em CR Acad. Sci. Paris Ser. B 275\/} (1972), 319--321.

\bibitem{deGennes}
{\sc De~Gennes, P.-G., and Prost, J.}
\newblock {\em The physics of liquid crystals}.
\newblock No.~83. Oxford university press, 1993.

\bibitem{dipasquale1}
{\sc Dipasquale, F., Millot, V., and Pisante, A.}
\newblock Torus-like solutions for the {Landau-de Gennes} model. {P}art {I}:
  the {L}yuksyutov regime.
\newblock {\em Archive for Rational Mechanics and Analysis 239\/} (2021),
  599--678.

\bibitem{dipasquale2}
{\sc Dipasquale, F.~L., Millot, V., and Pisante, A.}
\newblock Torus-like solutions for the {Landau-de Gennes model. Part II:
  Topology of S1}-equivariant minimizers.
\newblock {\em Journal of Functional Analysis 286}, 7 (2024), 110314.

\bibitem{dipasquale3}
{\sc Dipasquale, F.~L., Millot, V., and Pisante, A.}
\newblock Torus-like solutions for the {Landau-de Gennes model. Part III}:
  torus vs split minimizers.
\newblock {\em Calculus of Variations and Partial Differential Equations 63}, 5
  (2024), 136.

\bibitem{ericksen1991liquid}
{\sc Ericksen, J.~L.}
\newblock Liquid crystals with variable degree of orientation.
\newblock {\em Archive for Rational Mechanics and Analysis 113}, 2 (1991),
  97--120.

\bibitem{frank}
{\sc Frank, F.~C.}
\newblock On the theory of liquid crystals.
\newblock {\em Discussions of the Faraday Society 25\/} (1958), 19--28.

\bibitem{gartland2018scalings}
{\sc Gartland, E.~C.}
\newblock Scalings and limits of landau-de gennes models for liquid crystals: A
  comment on some recent analytical papers.
\newblock {\em Mathematical Modeling and Analysis 23}, 3 (2018), 414--432.

\bibitem{gm1999}
{\sc Gartland~Jr, E., and Mkaddem, S.}
\newblock Instability of radial hedgehog configurations in nematic liquid
  crystals under {Landau--de Gennes} free-energy models.
\newblock {\em Physical Review E 59}, 1 (1999), 563.

\bibitem{gz}
{\sc Geng, Z., and Zarnescu, A.}
\newblock Uniform profile near the point defect of {Landau-de Gennes} model.
\newblock {\em Calculus of Variations and Partial Differential Equations 62}, 1
  (2023), 3.

\bibitem{gm}
{\sc Golovaty, D., and Montero, J.~A.}
\newblock On minimizers of a {Landau--de Gennes} energy functional on planar
  domains.
\newblock {\em Archive for Rational Mechanics and Analysis 213}, 2 (2014),
  447--490.

\bibitem{han1994singular}
{\sc Han, Q.}
\newblock Singular sets of solutions to elliptic equations.
\newblock {\em Indiana University Mathematics Journal\/} (1994), 983--1002.

\bibitem{han1998geometric}
{\sc Han, Q., Hardt, R., and Lin, F.}
\newblock Geometric measure of singular sets of elliptic equations.
\newblock {\em Communications on Pure and Applied Mathematics 51}, 11-12
  (1998), 1425--1443.

\bibitem{hardt1999critical}
{\sc Hardt, R., Hoffmann-Ostenhof, M., Hoffmann-Ostenhof, T., and Nadirashvili,
  N.}
\newblock Critical sets of solutions to elliptic equations.
\newblock {\em Journal of Differential Geometry 51}, 2 (1999), 359--373.

\bibitem{hkl1}
{\sc Hardt, R., Kinderlehrer, D., and Lin, F.-H.}
\newblock Existence and partial regularity of static liquid crystal
  configurations.
\newblock {\em Communications in Mathematical Physics 105}, 4 (1986), 547--570.

\bibitem{hkl2}
{\sc Hardt, R., Kinderlehrer, D., and Lin, F.-H.}
\newblock Stable defects of minimizers of constrained variational principles.
\newblock {\em Annales de l'Institut Henri Poincar{\'e} C, Analyse non
  lin{\'e}aire 5}, 4 (1988), 297--322.

\bibitem{hardt1993harmonic}
{\sc Hardt, R., and Lin, F.~H.}
\newblock Harmonic maps into round cones and singularities of nematic liquid
  crystals.
\newblock {\em Mathematische Zeitschrift 213\/} (1993), 575--593.

\bibitem{hs}
{\sc Hardt, R., and Simon, L.}
\newblock Nodal sets for solutions of elliptic equations.
\newblock {\em Journal of Differential Geometry 30}, 2 (1989), 505--522.

\bibitem{hmp}
{\sc Henao, D., Majumdar, A., and Pisante, A.}
\newblock Uniaxial versus biaxial character of nematic equilibria in three
  dimensions.
\newblock {\em Calculus of Variations and Partial Differential Equations 56}, 2
  (2017), 1--22.

\bibitem{hqz}
{\sc Hu, Y., Qu, Y., and Zhang, P.}
\newblock On the disclination lines of nematic liquid crystals.
\newblock {\em Communications in Computational Physics 19}, 2 (2016), 354--379.

\bibitem{insz1}
{\sc Ignat, R., Nguyen, L., Slastikov, V., and Zarnescu, A.}
\newblock Uniqueness results for an {ODE} related to a generalized
  {Ginzburg--Landau} model for liquid crystals.
\newblock {\em SIAM Journal on Mathematical Analysis 46}, 5 (2014), 3390--3425.

\bibitem{insz2}
{\sc Ignat, R., Nguyen, L., Slastikov, V., and Zarnescu, A.}
\newblock Stability of the melting hedgehog in the {Landau--de Gennes} theory
  of nematic liquid crystals.
\newblock {\em Archive for Rational Mechanics and Analysis 215}, 2 (2015),
  633--673.

\bibitem{lamy2013}
{\sc Lamy, X.}
\newblock Some properties of the nematic radial hedgehog in the {Landau--de
  Gennes} theory.
\newblock {\em Journal of Mathematical Analysis and Applications 397}, 2
  (2013), 586--594.

\bibitem{lin1989nonlinear}
{\sc Lin, F.-H.}
\newblock Nonlinear theory of defects in nematic liquid crystals; phase
  transition and flow phenomena.
\newblock {\em Communications on Pure and Applied Mathematics 42}, 6 (1989),
  789--814.

\bibitem{lin1991nodal}
{\sc Lin, F.-H.}
\newblock Nodal sets of solutions of elliptic and parabolic equations.
\newblock {\em Communications on Pure and Applied Mathematics 44}, 3 (1991),
  287--308.

\bibitem{lin1991nematic}
{\sc Lin, F.~H.}
\newblock On nematic liquid crystals with variable degree of orientation.
\newblock {\em Communications on Pure and Applied Mathematics 44}, 4 (1991),
  453--468.

\bibitem{Lin-Liu}
{\sc Lin, F.-H., and Liu, C.}
\newblock Static and dynamic theories of liquid crystals.
\newblock {\em J. Partial Differential Equations 14}, 4 (2001), 289--330.

\bibitem{LW1}
{\sc Lin, F.-H., and Wang, C.-Y.}
\newblock Harmonic and quasi-harmonic spheres.
\newblock {\em Communications in Analysis and Geometry 7}, 2 (1999), 397--429.

\bibitem{LW2}
{\sc Lin, F.-H., and Wang, C.-Y.}
\newblock Harmonic and quasi-harmonic spheres, part ii.
\newblock {\em Communications in Analysis and Geometry 10}, 2 (2002), 341--375.

\bibitem{lyu}
{\sc Lyuksyutov, I.}
\newblock Topological instability of singularities at small distances.
\newblock {\em Zh. Eksp. Teor. Fiz 75\/} (1978), 358--360.

\bibitem{ma2012}
{\sc Majumdar, A.}
\newblock The radial-hedgehog solution in {Landau--de Gennes}' theory for
  nematic liquid crystals.
\newblock {\em European Journal of Applied Mathematics 23}, 1 (2012), 61--97.

\bibitem{mz}
{\sc Majumdar, A., and Zarnescu, A.}
\newblock {Landau--de Gennes} theory of nematic liquid crystals: the
  {Oseen--Frank} limit and beyond.
\newblock {\em Archive for Rational Mechanics and Analysis 196}, 1 (2010),
  227--280.

\bibitem{mg2000}
{\sc Mkaddem, S., and Gartland~Jr, E.}
\newblock Fine structure of defects in radial nematic droplets.
\newblock {\em Physical Review E 62}, 5 (2000), 6694.

\bibitem{nz}
{\sc Nguyen, L., and Zarnescu, A.}
\newblock Refined approximation for minimizers of a {Landau-de Gennes} energy
  functional.
\newblock {\em Calculus of Variations and Partial Differential Equations 47}, 1
  (2013), 383--432.

\bibitem{nomizu1973characteristic}
{\sc Nomizu, K.}
\newblock Characteristic roots and vectors of a diifferentiable family of
  symmetric matrices.
\newblock {\em Linear and Multilinear Algebra 1}, 2 (1973), 159--162.

\bibitem{palffy1994new}
{\sc Palffy-Muhoray, P., Gartland, E., and Kelly, J.}
\newblock A new configurational transition in inhomogeneous nematics.
\newblock {\em Liquid Crystals 16}, 4 (1994), 713--718.

\bibitem{su}
{\sc Schoen, R., and Uhlenbeck, K.}
\newblock A regularity theory for harmonic maps.
\newblock {\em Journal of Differential Geometry 17}, 2 (1982), 307--335.

\bibitem{su2}
{\sc Schoen, R., and Uhlenbeck, K.}
\newblock Regularity of minimizing harmonic maps into the sphere.
\newblock {\em Inventiones Mathematicae 78\/} (1984), 89--100.

\bibitem{ss}
{\sc Schopohl, N., and Sluckin, T.}
\newblock Hedgehog structure in nematic and magnetic systems.
\newblock {\em Journal de Physique 49}, 7 (1988), 1097--1101.

\bibitem{sonnet1995alignment}
{\sc Sonnet, A., Kilian, A., and Hess, S.}
\newblock Alignment tensor versus director: Description of defects in nematic
  liquid crystals.
\newblock {\em Physical Review E 52}, 1 (1995), 718.

\bibitem{ty}
{\sc Tai, H.-M., and Yu, Y.}
\newblock Pattern formation in {Landau--de Gennes} theory.
\newblock {\em Journal of Functional Analysis 285}, 1 (2023), 109923.

\bibitem{yu2020}
{\sc Yu, Y.}
\newblock Disclinations in limiting {Landau--de Gennes} theory.
\newblock {\em Archive for Rational Mechanics and Analysis 237}, 1 (2020),
  147--200.

\end{thebibliography}

\end{document}